# A Tight SDP Relaxation for MIQCQP Problems in Power Systems Based on Disjunctive Programming

Qifeng Li

*Abstract*—An optimization problem considering AC power flow constraints and integer decision variables can usually be posed as a mixed-integer quadratically constrained quadratic program (MIQCQP) problem. In this paper, first, a set of valid linear equalities are applied to strengthen the semidefinite program (SDP) relaxation of the MIQCQP problem without significantly increasing the problem dimension so that an enhanced mixed-integer SDP (MISDP) relaxation, which is a mixed-integer convex problem, is obtained. Then, the enhanced MISDP relaxation is reformulated as a disjunctive programming (DP) problem which is tighter than the former one, since the disjunctions are designed to capture the disjunctive nature of the terms in the rank-1 constraint about the integral variables. The DP relaxation is then equivalently converted back into a MISDP problem the feasible set of whose continuous relaxation is the convex hull of feasible region of the DP problem. Finally, globally optimal solution of the DP problem which is the tightest relaxation for the MIQCQP proposed in the paper is obtained by solving the resulting MISDP problem using a branch-and-bound (B&B) algorithm. Computational efficiency of the B&B algorithm is expected to be high since feasible set of the continuous relaxation of a MISDP sub-problem is the convex hull of that of the corresponding DP sub-problem. To further reduce the dimension of the resulting MISDP problem, a compact formulation of this problem is proposed considering the sparsity. An optimal placement problem of smart PV inverter in distribution systems integrated with high penetration of PV, which is an MIQCQP problem, is studied as an example. The proposed approach is tested on an IEEE distribution system. The results show that it can effectively improve the tightness and feasibility of the SDP relaxation.

## I. Nomenclature

| | |
|---|---|
| $c_s, c_s$ | Unit costs of smart inverters and conventional inverters respectively. |
| $P_{ij}, Q_{ij}$ | Active and reactive branch flow on branch $ij$. |
| $E, N$ | Branch set and node set of a certain power systems respectively. |
| $p_i^{Gen}, q_i^{Gen}$ | Active and reactive generations at bus $i$. |
| $p_i^{Load}, q_i^{Load}$ | Active and reactive load at bus $i$. |
| $p_i^{PV}$ | PV output at bus $i$. |
| $q_i^{Invt}$ | Reactive output of smart inverter at bus $i$. |
| $R$ | Rating of the substation transformer |
| $S_i^{Invt}, S_i^{PV}$ | Ratings of the smart inverter and PV panel at bus $i$ respectively. |
| $S_I$ | The set of integer variables. |
| $S_{LE}, S_{LI}$ | The sets of linear equalities and inequalities respectively. |
| $S_{QE}, S_{QI}$ | The sets of quadratic equalities and inequalities respectively. |
| $v_i$ | Square of magnitude of the voltage at bus $i$. |
| $\ell_{ij}$ | Square of magnitude of the current on branch $ij$. |
| $\alpha_i$ | Integer for $i \in S_I$. |

## II. Introduction

OPERATIONAL issues, e.g. overloading in conductors and power quality problems, in distribution systems integrated with high penetration of photovoltaic (PV) resources have been widely reported in literature [1] – [7]. Among these operational problems, the voltage violation is a severe one [3] - [7]. In 2009, a photovoltaic & storage integration research program conducted by EPRI identified common measures by which smart inverters may be integrated into utility systems [8]. The smart inverter volta/var control strategies for high penetration of PV on distribution systems were studied in [4]-[6]. The smart inverters are usually more expensive than the conventional ones due to the extra capability of reactive power support. The utilization of smart inverters in power systems raises an interesting issue: how to obtain the minimum investment of smart inverters to meet the volt/var control requirement.

The paper designs an optimization model for smart inverter placement minimizing the total investment of inverters (including smart inverters and regular inverters). The designed optimization model is a MIQCQP problem taking into account the AC power flow constraints. MIQCQP problems are hard to solve since they contain two kinds of non-convexities: integer variables and non-convex quadratic constraints (i.e. AC power flow equations) [9].

Recent years, numerous publications studied the convex relaxations of the AC power flow equations. This study is reflected in the research on convexification of the famous optimal power flow (OPF) problem (please refer to [10], [11] and the references therein). The semidefinite programming (SDP)

The author is with the Department of Electrical Engineering, Arizona State University, Tempe, AZ 85287 USA (e-mail: qifeng.li@asu.edu).

relaxation is one of the most popular convex relaxations used to convexify the OPF problem. Meanwhile, some researchers focused on convexification of the mixed-integer problems, e.g. security-constrained unit commitment [12], distribution systems reconfiguration [13], transmission system planning [14], reactive power planning [15] and optimal transmission switching [16]. Some of the references (e.g. [12]) treated the integral variables as continuous variables, which makes the mixed-integer problem have no evident difference from the continuous QCQPs. The others used branch-and-bound (B&B), branch-and-cut (B&C) algorithms or solvers with the former two algorithms implemented to obtain integral solutions for the integer variables.

The approach proposed in this paper is based on the SDP relaxation. There are a number of linear equalities in the optimization model for smart inverter placement problem. A cluster of linear equalities regarding the auxiliary variables are imposed on the relaxation to obtain a tighter SDP relaxation. Computational study shows that the linear equalities do not significantly increase the runtime of solving the problem.

A disjunctive nature of the rank-1 constraint $X = xx^T$ with respect to the integral variables is found in this paper. That is, for instance, $X_{ij} = x_i x_j$ implies $X_{ij} = x_i$ (if $x_j = 1$) and $X_{ij} = 0$ (if $x_j = 0$), where $x_j$ is a binary variable. This disjunctive property is hard to formulate in the conventional SDP framework. However, by properly designed disjunction terms, this nature can be perfectly captured in a disjunctive programming (DP) [17]-[20] framework. With this disjunctive property captured, the rank-1 constraint can be further approximated. As a result, the obtained DP problem is a tight relaxation for MIQCQPs, e.g. the optimal placement of smart inverters. Note that the disjunctions designed in [14] and [16] were used to deal with the on/off constraints instead of the disjunctive nature of the rank-1 constraint mentioned above. The order of on/off constraints is usually 1-higher than the order of other constraints in the studied problems. Disjunctive constraints are a widely used technique to reduce the order of the on/off constraints in literature [21] - [23].

The DP problem cannot be directly solved. A convex hull reformulation technique [17]-[20] is used to convert the DP relaxation back to a mixed-integer SDP problem. Feasible set of the continuous relaxation of the resulting MISDP problem is the convex hull of that of the DP problem. Note that, simply replace the integer variables with continuous variables, the resulting continuous problem is defined as the continuous relaxation of the original discrete problem in the paper. A B&B algorithm is used to achieve the global mix-integer solution for the DP problem by solving the resulting MISDP relaxation. At each node of the algorithm, feasible set of the continuous relaxation of a MISDP sub-problem is the convex hull of that of the corresponding DP sub-problem. This property may help increase the computational efficiency of the B&B algorithm.

Flow chart of the proposed convexification procedure for the optimal smart inverter placement problem studied in this paper is given in Fig. 1. In fact, this procedure can be applied to some MIQCQP cases that are more complex than the inverter placement problem. The extendibility of the proposed procedure is discussed in Section V. When the procedure is applied to a MIQCQP problem with general integral variables, a depth-first-based B&B algorithm is designed based on an observation in [20].

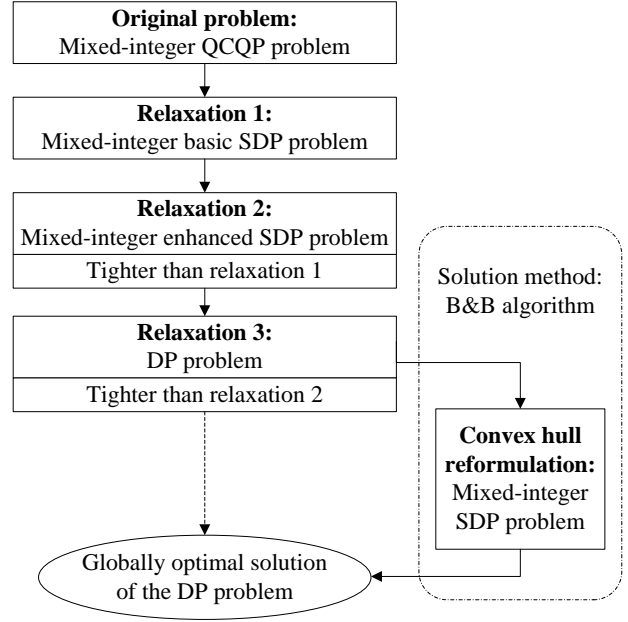

Fig. 1 Flow chart of the proposed convexification procedure for a MIQCQP problem.

The rest of the paper is organized as follow. In Section III, an optimal placement problem of smart PV inverters is modeled and its enhanced MISDP relaxation is proposed. In Section IV, DP reformulation of the enhanced MISDP relaxation is proposed and its convex hull is studied. Solving method and extendibility of the proposed approach are discussed in Section V and VI respectively. A case study is given in Section VII.

## III. MODELING OF THE OPTIMAL SMART PV INVERTER PLACEMENT PROBLEM AND ITS ENHANCED MISDP RELAXATION

### A. Mathematical Model

To determine the minimum total investment of inverters to meet the requirement of volt/var regulation, an optimization model is proposed as (1) where the AC power flow constraints are considered and the objective function as well as some constraints contain binary variables. As a matter of fact, a distribution system is usually operated in a radial typology. Hence, a branch flow model (BFM) [24] and [25] is used to describe the power flow in radial networks in the optimization model since the BFM is exact for a radial network and contains fewer non-convex quadratic constraints [7]. Although the application instance is in radial networks, the proposed method is still adaptable to a meshed network (please refer to the discussion in Section V). Note that there are mainly two BFMs and they were studied in [24] and [25] respectively. What used in the following optimization model is the one introduced in [24].

(MIQCQP) $\quad \min \ f = \sum_i \left[ c_s S_i^{Invt} + (1-\alpha_i) c_c S_i^{PV} \right]$ (1a)

s.t. $\quad p_i^{Gen} + p_i^{PV} - p_i^{Load} = \sum_{k:i \to k} P_{ik} - (P_{ji} - r_{ij} l_{ji})$ (1b)

$\quad q_i^{Gen} - q_i^{Load} + q_i^{Invt} = \sum_{k:i \to k} Q_{ik} - (Q_{ji} - x_{ij} l_{ji})$ (1c)

$\quad v_i = v_j + 2(r_{ij} P_{ij} + x_{ij} Q_{ij}) - (r_{ij}^2 + x_{ij}^2) \ell_{ij}$ (1d)

$\quad v_i \ell_{ij} = P_{ij}^2 + Q_{ij}^2$ (1e)

$\quad (q_i^{Invt})^2 + (p_i^{PV})^2 \leq (S_i^{Invt})^2 + (1-\alpha_i)(S_i^{PV})^2$ (1f)

$\quad -S_i^{Invt} \leq q_i^{Invt} \leq S_i^{Invt}$ (1g)

$\quad \alpha_i S_i^{PV} \leq S_i^{Invt} \leq \alpha_i M$ (1h)

$\quad 0 \leq \ell_{ij} \leq \overline{\ell}_{ij}$ (1i)

$\quad \underline{v}_i \leq v_i \leq \overline{v}_i$ (1j)

$\quad -0.6 R \leq p_{Gridt} \leq R$, (1k)

where the objective function represents minimizing the total investment of inverters in the system; (1b)-(1e) denote the branch flow constraints where (1b)-(1c) are linear and (1e) is quadratic; (1f) and (1g) represent the reactive capability constraints for the smart inverters which are quadratic inequalities and it is designed based on the research result of [6]; (1h) is the yes/no constraints for installing a smart inverter; (1i)-(1k) denote the system constraints. (1g) is dominated by (1f) when $\alpha_i = 1$ and conversely dominates (1f) when $\alpha_i = 0$. In many feeders, the only $p^{Gen}$ and $q^{Gen}$ are the active and reactive grid power.

Problem (1) is a MIQCQP problem. This is a relatively simple formulation for the optimal placement problem of smart PV inverters. Engineers may want to further consider some other elements, like capacitor banks and tap-changeable transformers as control variables [26], or solve the problem in a meshed distribution network. The proposed approach is in fact extendable to some more complex cases. Please refer to a discussion given in Section V.

### B. The Enhanced MISDP Relaxation

Let $x = [P, Q, p^{Gen}, q^{Gen}, v, \ell, S^{Invt}, q^{Invt}, \alpha]$ and replace the quadratic terms with auxiliary variables $X$ and omitting the rank-1 constraint $X = xx^T$, Problem (1) is relaxed to a basic MISDP problem [12] which is a mixed-integer convex problem (as shown in (2)). Note that the Shor's inequality [27] in (2g) is considered as part of the basic MISDP relaxation in the paper. Due to the existence of linear equalities (1b)-(1d), a set of linear equalities as shown in (3) are valid for strengthen the basic MISDP relaxation.

(MIESDP) $\quad \min \ f(x, X) = tr(Q_0 X) + c_0^T x$ (2a)

s.t. $\quad tr(Q_i X) + c_i^T x \leq b_i \quad (i \in S_{QI})$ (2b)

$\quad tr(Q_i X) + c_i^T x = b_i \quad (i \in S_{QE})$ (2c)

$\quad c_i^T x = b_i \quad (i \in S_{LE})$ (2d)

$\quad c_i^T x \leq b_i \quad (i \in S_{LI})$ (2e)

$\quad (\max\{0, \underline{x}_i\})^2 \leq X_{ii} \leq \max\{\underline{x}_i^2, \overline{x}_i^2\}$ (2f)

$\quad \begin{bmatrix} 1 & x^T \\ x & X \end{bmatrix} \succeq 0$ (2g)

$\quad c_i^T X c_j - (b_i c_j^T + b_j c_i^T) x + b_i b_j = 0$ (3a)

$\quad X c_i = b_i x \quad (i, j \in S_{LE})$ (3b)

$\quad c_i^T X c_j = b_i b_j$ . (3c)

where $\alpha_i$ in $x$ is integral for $i \in S_I$. Linear equalities (3) are respectively generated from the following relations implied by (2d) [28] and [29, Remark 13.4.1].

$$(c_i^T x - b_i)(c_j^T x - b_j) = 0$$
$$x c_i^T x = x b_i \quad (i, j \in S_{LE})$$
$$c_i^T x c_j^T x = b_i b_j$$

Note that, in the enhanced MISDP relaxation, which is named (MIESDP) in this paper, for the optimal placement problem of smart PV inverters, (3a) is preferred since it is tighter than (3b) and contains fewer constraints than (3c). With (2d), it is easy to verify that (3a) is equivalent to (3c). When $i = j$, both (3a) and (3c) are equivalent to (3b). Hence, (3a) and (3c) are tighter than (3b). When $i = j$, constraint (3a) becomes $c_i^T X c_i - 2 b_i c_i^T x + b_i^2 = 0$. Due to the Shor's inequality in (2g), $0 = c_i^T X c_i - 2 b_i c_i^T x + b_i^2 \geq c_i^T x x^T c_i - 2 b_i c_i^T x + b_i^2 = (c_i^T x - b_i)^2$, which means $c_i^T x - b_i = 0$. Thus, (2d) is redundant if (3a) is adopted in (MIESDP).

**Theorem 1.** The relaxation (MIESDP) is equivalent to the original problem (MIQCQP) if matrix $C$ which consists of the coefficient vector $c_i$ ($i \in S_{LE}$) of the linear equalities is a full-rank matrix.

Proof. $C$ is a full-rank matrix means $|S_{LE}| = n$. Suppose that constraint (3b) is imposed to generate the relaxation (MIESDP). $b_i x = x b_i = x c_i^T x = x x^T c_i$, as a result, $X c_i = x x^T c_i$ which is equivalent to $XC = xx^T C$. Both sides of $XC = xx^T C$ post-multiply by $C^{-1}$, then $X = xx^T$. If (3a) or (3c) is adopted instead of (3b), (MIESDP) is also exact since (3a) and (3c) are tighter than (3b). □

Note that, usually, the matrix $C$ for the optimal placement problem of smart PV inverters is not full-rank since $|S_{LE}| < n$. Generally, the higher the rank of $C$, the tighter the relaxation (MIESDP) is.

## IV. DP RELAXATION AND ITS CONVEX HULL

### A. Disjunctive Programming: An tighter Relaxation

To achieve a tighter relaxation of (MIQCQP), (MIESDP) is reformulated as a DP problem where the disjunctions and logic positions are designed based on the disjunctive nature introduced in Section I. The DP relaxation of (MIQCQP) is given as

(GDP) $\quad$ (2a) - (2c), (2e) – (2g), (3a) and

$$\underset{k \in D_i}{\vee} \begin{bmatrix} Y_{ik} \\ x_i = a_{ik} \\ X_{ij} = a_{ik} x_j \end{bmatrix} \quad (i \in S_I;\ j = 1, \ldots, n) \quad (4)$$

where $X_{ij}$, $x_i$ ($i, j = 1, \ldots, n$) are continuous variables, $Y_{ik} \in \{True, False\}$ and $\Omega(Y_i) = True$. If $x_g = \{0, 1, 2, 3\}$ for instance, the disjunction term (4) for $i = g$ becomes (5). The above DP problem is in its generalized form [19] and [20], therefore we named it (GDP).

$$\Omega(Y_g) = \vee \begin{bmatrix} Y_{g1} \\ x_g = 0 \\ X_{gj} = 0 \end{bmatrix} \vee \begin{bmatrix} Y_{g2} \\ x_g = 1 \\ X_{gj} = x_j \end{bmatrix} \vee \begin{bmatrix} Y_{g3} \\ x_g = 2 \\ X_{gj} = 2x_j \end{bmatrix} \vee \begin{bmatrix} Y_{g4} \\ x_g = 3 \\ X_{gj} = 3x_j \end{bmatrix}$$
$$(j = 1, \ldots, n) \quad (5)$$

As introduced in Section I, disjunctions (4) represent the entries of the rank-1 constraint $X = xx^T$ that is related to the integer variables. This relation is hard to directly formulate in the framework of the SDP relaxation. Therefore, (GDP) is tighter than (MIESDP). If the integral variables are non-binary, the effect of improving the tightness will be stronger since (4) will contain more constraints that is hard to model in (MIESDP). It is important to note that, in the above DP model, all the constraints inside and outside disjunctions (4) are convex.

*B. Convex Hull of (GDP)*

By adding some auxiliary variables, the disjunctions (4) can be equivalently reformulated as follow

$$\boldsymbol{y} = \sum_{k \in D_i} \boldsymbol{u}_{ik} \quad (i \in S_I) \quad (6a)$$

$$\lambda_{ik} \underline{\boldsymbol{y}} \leq \boldsymbol{u}_{ik} \leq \lambda_{ik} \overline{\boldsymbol{y}} \quad (k \in D_i;\ i \in S_I) \quad (6b)$$

$$\boldsymbol{A}_{ik} \boldsymbol{u}_{ik} = \lambda_{ik} \boldsymbol{B}_{ik} \quad (k \in D_i;\ i \in S_I) \quad (6c)$$

$$\sum_{k \in D_i} \lambda_{ik} = 1 \quad (i \in S_I) \quad (6d)$$

where $\boldsymbol{y} = [X_{11}, X_{12}, \ldots, X_{n,n-1}, X_{n,n}, x_1, \ldots, x_n]'$ is an intermediate variable vector, $\overline{\boldsymbol{y}}$ and $\underline{\boldsymbol{y}}$ are upper and lower bounds of $\boldsymbol{y}$ respectively; $\boldsymbol{u}_{ik}$ is the vector of auxiliary variables and $\lambda_{ik} = \{0, 1\}$ ($k \in D_i;\ i \in S_I$). (6c) denotes the relations within the $k$th term of the $i$th disjunction which is described in (4). For more details, please refer to the appendix section.

Equations (6) define a discrete feasible set in the space ($\boldsymbol{x}$, $\boldsymbol{X}$, $\boldsymbol{u}$, $\boldsymbol{\lambda}$) whose projection onto the ($\boldsymbol{x}$, $\boldsymbol{X}$)-space is exactly the feasible region determined by disjunctions (4). When a continuous relaxation is applied to (6) (i.e., "$\lambda_{ik} = \{0, 1\}$" is replaced with "$0 \leq \lambda_{ik} \leq 1$"), the feasible set defined by it becomes a convex set in the space ($\boldsymbol{x}$, $\boldsymbol{X}$, $\boldsymbol{u}$, $\boldsymbol{\lambda}$) since all equations in (6) are linear. Dramatically, the projection of this convex set onto the ($\boldsymbol{x}$, $\boldsymbol{X}$)-space is the convex hull of the feasible region of (4). The proof of the above statement which can be regarded as an extension of the related works in [17]-[20].

**Remark.** The projection of the feasible region of the continuous relaxation of (6) onto the ($\boldsymbol{x}$, $\boldsymbol{X}$)-space is the convex hull of that defined by (4). It does not necessarily mean that the convex feasible region specified by equations (6) with "$0 \leq \lambda_{ik} \leq 1$" is the convex hull of that in ($\boldsymbol{x}$, $\boldsymbol{X}$, $\boldsymbol{u}$, $\boldsymbol{\lambda}$)-space defined by equations (6) with "$\lambda_{ik} = \{0, 1\}$". In fact, it is most probably not. Just like, simply making $\alpha_i$ ($i \in S_I$) in (MIESDP) continuous will not result in the convex hull of its feasible set.

The optimization problem consists of (2a) - (2c), (2e) – (2g), (3a) and (6) with "$\lambda_{ik} = \{0, 1\}$" is named as (CH-MIESDP) while that consists of (2a) - (2c), (2e) – (2g), (3a) and (6) with "$0 \leq \lambda_{ik} \leq 1$" is referred to as (CH-ESDP). In the ($\boldsymbol{x}$, $\boldsymbol{X}$)-space, the feasible set of (CH-ESDP) is the convex hull of that of (CH-MIESDP) as well as that of (GDP) ((CH-MIESDP) is equivalent to (GDP)).

*C. An Illustrative Example*

To intuitively reveal the importance of obtaining the convex hull in improving efficiency of the B&B algorithm, an illustrative example is provided in this subsection. Consider the following mixed-integer quadratic inequality

$$x^2 + 2xy + y^2 \leq 1 \quad (7)$$

where $y = \{0, 1\}$. The feasible set of this mixed-integer inequality is the point (0, 1) and the segment of $x$-axis where $|x| \leq 1$ as shown in Fig. 2. It is easy to show that feasible set of the continuous relaxation is the semicircular region while the convex hull is denoted as the shaded triangular area.

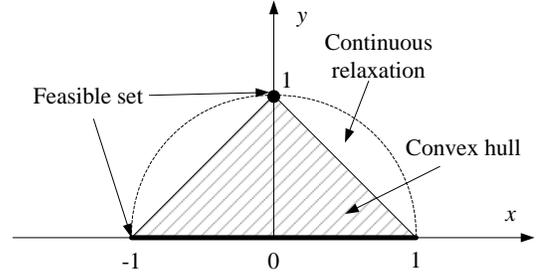

Fig. 2. The illustrative example.

Suppose that problem (7) is the sub-problem at certain node of the B&B algorithm, the computer searches for the optimal solution in the semicircle area. When the convex hull reformulation of the inequality is solved using the B&B algorithm, the searching area is the triangular region which is smaller than the semicircular region. It can be expected that the computational efficiency will be improved significantly with a much smaller searching area [20].

The efficiency of the convex hull reformulation can also be understood through another perspective. Valid cutting planes are usually added to improve the computational efficiency of the B&B algorithm by shrinking the searching region (i.e. the well-known B&C algorithm). However, there is no cutting plane available for enhancing the B&B algorithm when (CH-MIESDP) is solved, since any cutting plane that cuts off a part of the convex hull will also cut off some part of the (GDP) feasible set.

## D. Comparison of the Tightness

Tightness of the convex relaxations for the optimal placement of smart PV inverters is compared in Fig. 3. (GDP) is the tightest relaxation of the original problem proposed in this paper, however it is not convex. The objective of the paper is to obtain the global solution of (GDP) by solving its equivalent problem (CH-MIESDP) using a B&B algorithm. A basic B&B algorithm solves the continuous relaxation of a discrete subproblem at each node of the algorithm. The algorithm is expected to be effective since the continuous relaxation of (CH-MIESDP) offers the convex hull of (GDP) as discussed in Subsection III-B.

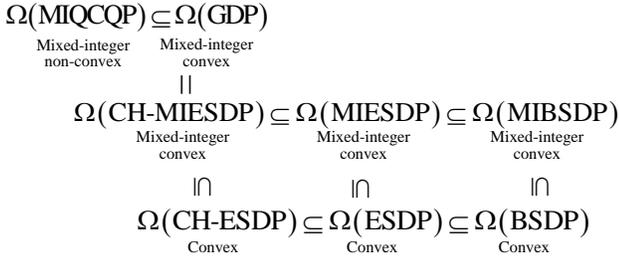

Fig. 3. Comparison of the tightness. $\Omega(\cdot)$ denotes the feasible set of a given problem. (MIBSDP) is the problem described by (2). (ESDP) and (BSDP) are the continuous relaxations of (MIESDP) and (MIBSDP) respectively. (CH-SDP) represents the SDP relaxation whose feasible region is the convex hull of that of (GDP).

However, the convex hull formulation increases the number of variables by $(k*m*(t+1))$ where $t$ is the total number of variables in the SDP relaxation, $m = |S_I|$ and $k = |D_i|$. If the integer variables are binary, then $k = 2$. If both of $t$ and $i$ are large values, the increase in problem size may weaken the benefit of improving the computational efficiency. To overcome this drawback, a compact formulation of (6) is proposed in the next subsection.

## E. A Compact Formulation Considering Sparsity

The matrix equalities in (6c) are the reformulations of the corresponding simple equality constraints in (4). Hence, $A_{ik}$ in (6c) is a highly sparse matrix where a large number of columns are 0 vectors. As a result, equations (6) only influence the variables that appear in (4). Therefore, $y$ can only represents the variables that appear in (4) so that the zero columns in $A_{ik}$ can be removed. As a result, the compact formulation of (6) is given in (8) which is similar to (6). However, the dimensions of $y$, $u$ and $A$ therein are much lower than those in (6).

$$\tilde{y} = \sum_{k \in D_i} \tilde{u}_{ik} \quad (i \in S_I) \tag{8a}$$

$$\lambda_{ik} \underline{\tilde{y}} \leq \tilde{u}_{ik} \leq \lambda_{ik} \overline{\tilde{y}} \quad (k \in D_i; i \in S_I) \tag{8b}$$

$$\tilde{A}_{ik} \tilde{u}_{ik} = \lambda_{ik} \tilde{B}_{ik} \quad (k \in D_i; i \in S_I) \tag{8c}$$

$$\sum_{k \in D_i} \lambda_{ik} = 1 \ (i \in S_I) \tag{8d}$$

where $\tilde{y} = [X_{i,1}, X_{i,2}, \ldots, X_{i,n}, x_1, \ldots, x_n]'$ ($i \in S_I$), $\overline{\tilde{y}}$ and $\underline{\tilde{y}}$ are upper and lower bounds of $\tilde{y}$ respectively. The problem consists of (2a) - (2c), (2e) – (2g), (3a) and (8) is denoted as (CH-MIESDP) in the paper.

**Theorem 2**. (CH-MIESDP) is equivalent to the optimization problem consists of (2a) - (2c), (2e) – (2g), (3a) and (6).

Proof. From the definitions of $y$ and $\tilde{y}$, it can be shown that

$$y = \begin{bmatrix} \mathbf{0}_{1 \times \left(\frac{(n-m)(3n-m+2)}{2}\right)} & \tilde{y}' \end{bmatrix}'$$

$$A_{ik} = \begin{bmatrix} \mathbf{0}_{\left(\frac{(n-m)(3n-m+2)}{2}\right) \times \left(\frac{(n-m)(3n-m+2)}{2}\right)} & \mathbf{0}_{\left(\frac{(n-m)(3n-m+2)}{2}\right) \times \left((n+1)m - \frac{(n-m)^2}{2}\right)} \\ \mathbf{0}_{\left((n+1)m - \frac{(n-m)^2}{2}\right) \times \left(\frac{(n-m)(3n-m+2)}{2}\right)} & \tilde{A}_{ik} \end{bmatrix}$$

$$B_{ik} = \begin{bmatrix} \mathbf{0}_{1 \times \left(\frac{(n-m)(3n-m+2)}{2}\right)} & \tilde{B}_{ik}' \end{bmatrix}'$$

where $\mathbf{0}_{\cdot \times \cdot}$ is a $(\cdot \times \cdot)$-dimension zero matrix. As a result, it suffices to verify that a point in the $(x, X)$-space that satisfies (6) will satisfy (8), and vise versa. □

As a result, number of the auxiliary variables of the problem reduces from $(k*m*(t+1))$ to $(k*m*((n-m/2)*m+n+1))$ if the compact formulation (8) is used, where $t = (n+1)*n$ and $n$ is the number of variables in $x$. When $n \gg m$, the compactness of (8) is high.

## V. SOLVING METHODOLOGY

### A. B&B Algorithm for the Binary Case

When a B&B algorithm is used to solve an integer programming problem, some settings may affect the computational performance drastically [30]. These settings include node selecting strategy (i.e. depth-first search plus backtracking and breadth-first search) and strategy for branching variable selection (i.e. choosing the next integer variable on which to branching) [31]. Users can choose these strategies based on the problem they need to solve. However, there is no universal rule for making the choices. Some may select the branching variable with the lowest or highest objective value [30] while others may choose the smallest or largest fraction value [32].

In the case study section of this paper, a standard B&B algorithm provided by a built-in solver of MATLAB, BNB, is used. With BNB, one can choose the node selecting strategy expediently. For further information about BNB, please refer to the help text of MATLAB and [33]. At each node, a SDP solver is called to obtain the bounds for the corresponding subproblem.

### B. A Branching Strategy for the General Integer Case

If integer variables of the MIQCQP problem are non-binary, for instance capacitor banks and tap-changeable transformers are treated as controllable reactive resources [26] in the smart inverter placement problem, then $D_i$ in (8d) may contain more than 2 terms. In its relaxation (CH-MIESDP), a binary variable $\lambda_{ik}$ is used to represent one term in a disjunction, which means the integer variables in (CH-MIESDP) become binary. How-ever, this does not imply that the general integer case is completely reduced to a binary case. To achieve better computational performance, one may need to come up with strategies for both of choosing branching variable and branching term for each general integer variable.

After solving the continuous relaxation of (CH-MIESDP), the obtained value for $\lambda_{ik}$ denote the "closeness" between the obtained optimal solution and the related disjunction term [20]. Inspired by the above observation, a depth-first-based B&B algorithm for solving (CH-MIESDP) is designed in the paper which offers priority to search for the optimal solution in the disjunction term that is closest to the optimal solution obtained at the previous node (parent node). The main steps of the designed B&B algorithm are given in Fig. 4.

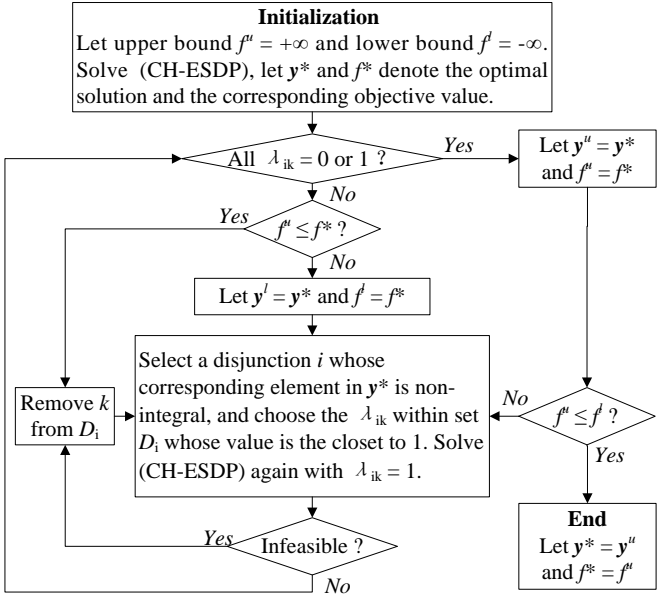

Fig. 4. Flow chart of the proposed B&B algorithm. $y = (\lambda_{mk}, \ldots, \lambda_{nk}, x)$, where $x$ denotes the vector of continuous variables and $k \in D_i$ as defined in (6). The compact formulation of (CH-ESDP) in (8) is recommended.

## VI. DISCUSSION: EXTENDIBILITY TO MORE COMPLEX CASES

Theoretically, the proposed approach is capable for some MIQCQP cases that are more complex than the smart inverter placement problem described in Section III. They are

**Similar cases in meshed networks:**

In a meshed network, the bus injection model (BIM) [10], namely the conventional formulation in rectangular coordinates, is usually used to describe the power flow equations, since the BFM is not valid to meshed networks. The power flow formulation in BIM is still a set of quadratic equalities. The SDP relaxation for the power flow in BIM can also be expressed as (2c). That means formulation (2) is also capable for a MIQCQP problem in meshed networks.

In an optimization model considering AC power flow constraints which are formulated using BIM, there may not be linear equality constraints. As a results, the linear equalities in (3) may not be valid for this case. However, the upper bound constraints for bus voltage magnitudes in the BIM, i.e. $\mathrm{Re}[V]^2 + \mathrm{Im}[V]^2 \leq \bar{V}^2$, are convex quadratic constraints. Some semidefinite inequalities that stem from these convex quadratic constraints may be valid for strengthening the SDP relaxation in the non-iterative framework [34].

**Cases where the quadratic terms contain integer variables:**

In problem (1), there is no quadratic term that contains integer variables. However, some MIQCQP problems in power systems contain integer variables in the quadratic constraints, like security-constrained unit commitment [12]. When defining the problem in Subsection III-B, the vector $x$ is defined including the integer variables, which means the potential quadratic terms that contain integer variables are replaced by the related entries in $X$. Hence, it suffices to show that the proposed method is valid for cases where the quadratic constraints contain integer variables.

**Cases with general integer variables:**

As introduced in Subsection V-B, the general integer variables means the non-binary integer variables. The disjunction instance given in (5) perfectly shows the capability of the proposed method for such cases. For better solving these cases, a potential B&B algorithm is designed in Subsection V-B.

**Cases for which the iterative constraints are effective:**

To mitigate the inexactness issue of the basic form of the SDP relaxation, the methods of utilizing some iteratively generated valid inequalities [35] to obtain tight SDP relaxations attract some researchers' attention in power systems. As reported by these references, the iterative methods work quite well in solving some unit commitment cases [36] as well as certain OPF cases [37].

In fact, the proposed approach offers an iterative computing architecture, which means it is possible to incorporate some valid inequalities into the proposed algorithm to achieve even tighter SDP relaxations. As stated in Subsection IV-C, the fact that feasible set of (CH-ESDP) is the convex hull of that of (GDP) means there is no valid cutting plane for obtaining the global solution of (GDP). However, (GDP) itself is a relaxation of (MIQCQP). It is still valuable to explore valid cutting planes that can effectively cut off the feasible set of (GDP) without harming that of (MIQCQP).

## VII. CASE STUDY

### A. Solving method

Solvers that can commendably solve a mixed-integer semidefinite programming problem are not currently available. BNB is an implementation of a standard branch & bound algorithm for mixed-integer convex programming and relies on external solvers for solving the node problems. Thus, in this case study, the B&B framework provided by BNB is used and MOSEK [38] is called for obtaining the lower bound of each node problem which is an SDP problem. The above solving procedure is implemented in YALMIP [39], a MATLAB optimization toolbox.

## B. The Test-bed system

The proposed approach is tested on the IEEE 13-bus standard feeder of which the typology is given in Fig. 5. Five commercial PV systems are assumed to connect to the feeder. The locations and power ratings of these PVs are also shown in Fig. 5 and Table I respectively. Total active load is 3.266 MW, then the penetration of PV is 73.5%. It is assumed that $c_s = 1.5c_i$.

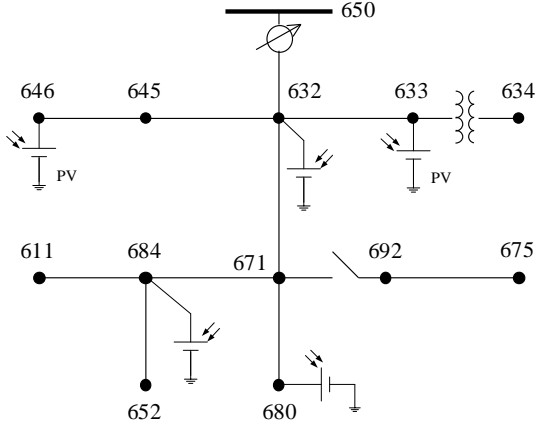

Fig. 5. IEEE 13-bus feeder.

TABLE II
LOCATIONS AND RATINGS OF PV SYSTEMS

| Node # | PV Rating (kW) |
|---|---|
| 632 | 200 |
| 633 | 600 |
| 646 | 400 |
| 680 | 700 |
| 684 | 500 |

## C. Results and Analysis

Optimal objective value (OOV, is called lower bound in some references) of each relaxation is shown for comparing the tightness while the ranks as well as the maximum entry of the error matrices ($X - xx^T$) are reported for comparing the feasibility of the solutions. The runtimes in per unit are compared. Since the runtime depends on not only the computer configurations but also coding skills, it is more clear and makes more sense to report the per unit runtimes using runtime of the basic case as the base. In this case study, the MISDP problem described in (2), namely (MIBSDP), is used as the basic case. Results are tabulated in Table II.

TABLE II
RESULTS OF CASE STUDY

| Relaxation | OOV | Max. Entry of Error Matrix | Rank of Error Matrix | Solver Time in p.u. |
|---|---|---|---|---|
| MIBSDP | 2.4 | 3.98 | 66 | 1 |
| MIESDP | 2.4 | 2.5639 | 66 | 0.96504 |
| CH-MIESDP | 3.6 | 2.4767 | 62 | 1.048317 |

It can be observed from Table II that both of the valid linear equalities and disjunctive formulation can improve the tightness of the SDP relaxation and feasibility of the obtained solution for the smart PV inverter placement problem which is a MIQCQP problem. The convex hull reformulation on one hand increases the dimension of the problem, tightness of the convex hull heightens the computational efficiency on the other hand. The computational result shows that the runtime of solving the convex hull reformulation taking into account the linear equalities (3a) is almost the same as that of solving the basic mixed-integer SDP relaxation of the inverter placement problem.

## VIII. CONCLUSION AND FUTURE WORK

The paper proposes a tight SDP relaxation for MIQCQP problems in power systems based on disjunctive programming where the disjunctive nature in the rank-1 constraints can be properly captured. A set of linear equalities is also used to tighten the relaxation. The proposed method is applied to obtaining the optimal placement of smart PV inverters in distribution systems integrated with high penetration of PV. The case study shows that the proposed approach can improve both tightness and feasibility of the SDP relaxation for the studied problem without remarkably increasing the computational burden.

In future work, the proposed approach will be applied to meshed distribution systems as well as some other MIQCQPs in transmission systems like SCUC. Efficiency of the method will also be tested on larger systems.

## IX. APPENDIX: PROOF OF (6)

Although the convex hull reformulation of the disjunctive programming has been well studied in literature [17]-[20], it is still necessary to guarantee that the application of this reformulation in the SDP relaxation of problems in power systems is correct. As a result, the following proof is provided:

Since the constraints in disjunctions (4) are linear, each term of a disjunction is a convex subset. The convex hull of one disjunction in (4) can be expressed as a convex combination of these convex subsets where all multipliers $\lambda_k$ are nonnegative and sum to 1. Hence, for $i \in S_I$, the convex combination of the subsets is

$$\lambda_{ik}\underline{y} \leq \lambda_{ik}y \leq \lambda_{ik}\overline{y} \quad (k \in D_i)$$
$$\lambda_{ik}A_{ik}y = \lambda_{ik}B_{ik} \quad (k \in D_i)$$
$$\sum_{k \in D_i}\lambda_{ik} = 1$$

where $A_{ik}y = B_{ik}$ is exactly the matrix form of the linear equalities within the disjunctions in (4), therefore $A_{ik}$ is highly sparse.

Replace the quadratic term $\lambda_{ik}y$ with auxiliary variable $u_{ik}$ for $k \in D_i$, the following relations is obtained.

$$u_{ik} = \lambda_{ik}y \quad (A.1)$$
$$\lambda_{ik}\underline{y} \leq u_{ik} \leq \lambda_{ik}\overline{y} \quad (k \in D_i) \quad (A.2)$$

$$A_{ik} u_{ik} = \lambda_{ik} B_{ik} \quad (k \in D_i) \quad (A.3)$$

$$\sum_{k \in D_i} \lambda_{ik} = 1 \quad (A.4)$$

Due to (A.4), it is straightforward to show that (A.1) is equivalent to

$$y = \sum_{k \in D_i} \lambda_{ik} y = \sum_{k \in D_i} u_{ik}. \quad (A.5)$$